\documentclass[a4paper,11pt,centertags]{amsart}
\usepackage{amsmath}
\usepackage{amsfonts}
\usepackage{amssymb}
\usepackage[latin1]{inputenc}
\usepackage{amsthm}
\usepackage[arrow, matrix, curve]{xy}

\theoremstyle{plain}
\newtheorem{theorem}{Theorem}[section]
\newtheorem{proposition}[theorem]{Proposition}
\newtheorem{lemma}[theorem]{Lemma}
\newtheorem{corollary}[theorem]{Corollary}

\theoremstyle{definition}

\newtheorem{convention}[theorem]{Convention}
\newtheorem{definition}[theorem]{Definition}

\theoremstyle{remark}

\newtheorem{example}[theorem]{Example}

\begin{document}

\setlength{\parindent}{0cm} \setcounter{tocdepth}{2}
\newcommand{\R}{\mathbb R}
\newcommand{\C}{\mathbb C}
\newcommand{\Z}{\mathbb Z}
\newcommand{\Q}{\mathbb Q}
\newcommand{\tr}{{\rm tr}}
\newcommand{\G}{\Gamma}
\newcommand{\K}{\mathbb K}
\newcommand{\M}{\mathcal M}
\newcommand{\X}{\mathcal X}
\renewcommand{\index}{\textbf}
\renewcommand{\o}{\omega}
\renewcommand{\L}{\mathfrak}
\newcommand{\Lg}{\widehat{L(\Li g)}}
\renewcommand{\d}{\partial}
\newcommand{\ra}{\rightarrow}
\newcommand{\xra}{\xrightarrow}
\newcommand{\+}{\oplus}
\newcommand{\x}{\otimes}
\renewcommand{\l}{\langle}
\renewcommand{\phi}{\varphi}
\renewcommand{\r}{\rangle}
\newcommand{\V}{\mathcal V}
\newcommand{\tri}{\triangle}
\newcommand{\sx}{| \negthickspace\negthickspace \times}
\def\cprime{$'$}

\bibliographystyle{plain}

\title[Quasimorphisms]{Reconstructing Quasimorphisms From Associated Partial Orders and A Question Of Polterovich}
\author{Gabi Ben Simon and Tobias Hartnick}

\begin{abstract} We show that every continuous homogeneous quasimorphism on a finite-dimensional $1$-connected simple Lie group arises as the relative growth of  some continuous bi-invariant partial order on that group. More generally we show, that an arbitrary homogeneous quasimorphism can be reconstructed as the relative growth of a partial order subject to a certain sandwich condition. This provides a link between invariant orders and bounded cohomology and allows the concrete computation of relative growth for finite dimensional simple Lie groups as well as certain infinite-dimensional Lie groups arising from symplectic geometry.
\end{abstract}
\maketitle

\section{Introduction}

In this article we observe a new relation between two different well-known structures on Lie groups. The one side of our correspondence is formed by continuous invariant partial orders. Here a partial order $\leq$ on a topological group $G$ is called \emph{invariant} (or bi-invariant), if for all $g,h,k \in G$ the relation $g \leq h$ implies both $kg \leq kh$ and $gk \leq hk$. This means that the associated \emph{order semigroup}
\[G^+ := \{g\in G \,|\, g \geq e\}\]
is a conjugation-invariant pointed (i.e. $G^+ \cap (G^+)^{-1} = \{e\}$) monoid. Then $\leq$ is called \emph{continuous} if $G^+$ is closed in $G$ and locally topologically generated (i.e. for every identity neighbourhood $U$ in $G$ the intersection $U \cap G^+$ generates a dense subsemigroup of $G^+$). Such orders will be related to continuous homogeneous quasimorphism, i.e. continuous maps $f: G \to \R$ satisying $f(g^n) = nf(g)$ for all $g \in G$ and $n \in \Z$, for which the function $f(gh) - f(g)-f(h)$ is bounded on $G^2$.  Both sides of the correspondence individually are well-studied; for finite-dimensional simple Lie groups there are classifications of both (see \cite{Gichev} for invariant orderings and \cite{BuMo} for quasimorphisms). An immediate consequence of these classification results is the following proposition:
\begin{proposition}\label{PropMotivation}
For a finite-dimensional $1$-connected simple Lie group $G$ the following are equivalent:
\begin{itemize}
\item[(i)] There exists a non-trivial continuous invariant partial order $\leq$ on $G$.
\item[(ii)] There exists a non-zero continuous homogeneous quasimorphism on $G$.
\end{itemize}
\end{proposition}
The main result of this article states that we can actually use the continuous invariant partial orders to construct the corresponding continuous homogeneous quasimorphisms explictly. For this we use the machinery of relative growth as introduced in \cite{EP}: Given any invariant partial order $\leq$ on a group $G$ we define the associated set of \emph{dominants} in $G$ to be
\[G^{++} = \{g \in G^+\setminus\{e\}\,|\, \forall h \in G \exists n \in \mathbb N: g^n \geq h\}.\]
We call an invariant order \emph{admissible} if $G^{++} \neq \emptyset$. For a fixed dominant element $g \in G^{++}$ we define the \emph{relative growth}
\[\gamma(g, \cdot): G \to \R\]
by
\[\gamma(g,h) = \lim_{n \to \infty} \frac{\min \{p \in \Z\,|\, g^p \geq h^n\}}{n}.\]
Then we provide the following explicit correspondence:
\begin{theorem}\label{IntroMain}
Let $G$ be a $1$-connected simple Lie group. Then any  continuous invariant order $\leq$ on $G$ is admissible, and the relative growth of any such order is given by
\[\gamma(g,h) = \pm \frac{f(h)}{f(g)} \quad (g \in G^{++}, h \in G),\]
where $f: G \to \R$ is any generator of the one-dimensional space of homogeneous quasimorphisms on $G$. 
\end{theorem}
To the best of our knowledge this is the first results which provides a correspondence between invariant order structures and quasimorphisms. In modern language, the two sides of the correspondence are given by Lie semigroups \cite{HHL, Neeb} and continuous bounded cohomology classes \cite{BC, BuMo}, respectively. In fact, the study of continuous invariant orders reduces to the more classical subject of invariant cones in Lie algebras. Interest in such invariant cones first arose in the context of infinite-dimensional representation theory and mathematical physics (in particular, general relativity)  \cite{Vinberg, Olshanski, Paneitz}. By now Lie semigroups have found applications in areas as diverse as logic and geometric control theory (see \cite{Hofmann} for a historical overview). On the other hand, bounded cohomology in general and quasimorphisms in particular are an indispensable tool in modern geometric group theory. Some articles of particular relevance to the present work are \cite{BargeGhys, Bavard, BuMo, Surface}. We hope that the present work will initiate more interaction between these two rich and traditional areas of topological group theory. We would like to point out that the theory of relative growth was originally developed in  \cite{EP} in a completely different context, namely the study of infinite-dimensional Lie groups arising from problems in contact and symplectic geometry.\\

In order to motivate the first step in the proof of Theorem \ref{IntroMain}, we consider a purely algebraic variant of that theorem, which applies to general groups $G$ and arbitrary
homogeneous quasimorphism $f: G \to \R$. For such a pair $(G,f)$ one can always construct an invariant partial order $\leq_f$ on $G$ by demanding that 
\[g <_f h \Leftrightarrow f(h^{-1}g) < - D(f),\]
where
\[D(f):= \sup_{g,h \in G} (f(gh) -f(g) -f(h))\]
denotes the \emph{defect} of $f$. Using the theory of relative growth we can show that $\leq_f$ actually determines $f$ up to a positive multiplicative constant. Indeed, we have:
\begin{proposition}\label{dull}
Let $G$ be a  group, $f: G \to \R$ a homogeneous quasimorphism and $\leq_f$ as above. Then $\leq_f$ is admissible, and for any  $g
\in G^{++}$ and $h \in G$ the corresponding relative growth is given by
\[\gamma(g,h) = \frac{f(h)}{f(g)}.\]
In particular, up to positive multiple every quasimorphism arises as the relative growth of some partial order with respect to any dominant.
\end{proposition}
Proposition \ref{dull} is an interesting observation in its own right, but it does \emph{not} provide continuous orders. In fact, the order semigroup of $\leq_f$ will not even be connected. We will thus need a stronger version of Proposition \ref{dull}, which can be used to compute the relative growth of continuous orders as well. Namely, we will show in Proposition \ref{Sandwich} below that a homogeneous quasimorphism $f$ can be recovered as relative growth from \emph{any} invariant partial order, which agrees with $\leq_f$ up to some bounded error.  In a second step we have to obtain explicit descriptions of all continuous quasimorphisms and continuous invariant orders on $1$-connected simple Lie groups. The main work then lies in the third and final step, where we use the results of Step 2 in order to verify that the continuous quasimorphisms are related to the continuous orders in such a way that Proposition \ref{Sandwich} applies. While the first step uses only elementary methods, the other two steps depend on an in depth understanding of the fine structure of the Lie groups under consideration.\\

This article is organized as follows: In Section \ref{SecHermitian} we recall the structure of those $1$-connected simple Lie groups, which admit continuous invariant partial orders. These turn out to be Hermitian. We obtain explicit descriptions of both continuous homogeneous quasimorphisms and continuous invariant partial orders on such groups. In Section \ref{SecMain} we prove Theorem \ref{IntroMain} along the lines explained above. We first provide the necessary generalization of Proposition \ref{dull} and use it to reduce the statement of the theorem to an estimate on the values of the quasimorphism in question. This estimate willl be established separately for the contribution coming from a maximal compact subgroup and a complementary non-compact contribution. In a final subsection we indicate briefly how to generalize our results beyond the simply-connected case. The concluding Section \ref{SecApps} discusses various applications and extensions of the main result. Following \cite{EP} we introduce the notion of an order space, which is a certain metric space associated to an ordered group. We explain how our main results allow one to compute the order space of $1$-connected simple Hermitian Lie group for suitable orderings. This answers in particular a question of Polterovich, which was the starting point for the present article. We then discuss possible extensions of our results to infinite-dimensional Lie groups arising in symplectic geometry. Again, we are able to compute certain order spaces, and the results are in strong contrast to existing results about similar infinite-dimensional groups in the symplectic context.\\

\begin{convention} In order to avoid tedious repetitions,  throughout the body of this article all homogeneous quasimorphisms are assumed to be continuous. (Note that in fact any  homogeneous quasimorphism on a finite-dimensional simple Lie group is automatically continuous  \cite{Shtern}.) 
\end{convention}

{\bf Acknowledgement:} We cordially thank Leonid Polterovich for suggesting the problem of computing the order spaces of simple Lie groups, which was the starting point for
this paper, and for pointing out the applications of our criterion to groups of Hamiltonian diffeomorphisms. We also thank Marc Burger for a number of
useful discussions concerning Hermitian Lie groups and their quasimorphisms. This article would not have been
possible without the competent guidance of Karl Heinrich Hofmann through the vast literature on Lie semigroups, which is gratefully acknowledged. The second-named
author was partially supported by Swiss National Science Foundation (SNF), grant PP002-102765.

\section{Quasimorphisms and partial orders on Hermitian Lie groups}\label{SecHermitian}

\subsection{The structure of Hermitian Lie groups}
Let $G$ be a $1$-connected simple real Lie group and  $G_0 := {\rm Ad}(G)$ so that $\pi: G \ra G_0$ is a universal covering map. Fix a maximal compact subgroup $K_0\subset G_0$ and define $K := \pi^{-1}(K_0)$. Then $\X := G/K = G_0/K_0$ is a symmetric space and $G$ is called Hermitian if the space $\Omega^2(\X)^G$ of $G$-invariant $2$-forms on $\X$ is non-trivial. In this case, actually, $\Omega^2(\X)^G  \cong  \R$. It was proved by Vinberg \cite{Vinberg} that among simple Lie groups only the Hermitian ones can admit continuous invariant partial orders. We will thus focus on such Lie groups in the sequel. In order to fix our notation we briefly recall the structure theory of $1$-connected simple Hermitian Lie groups. For more details the reader is asked to consult \cite{KW}, \cite[Chapter III]{Koranyi} and (regarding compact
Lie groups) \cite[Chapter IV]{Knapp}. Throughout, the Lie algebra of a Lie group is denoted by the corresponding small gothic letter; a subscript $\C$ indicates complexification.

\begin{itemize}
\item $\L k$ decomposes as $\L k = \L{z(k)} \oplus \L k'$, where  $\L{z(k)}$ denotes the $1$-dimensional center of $\L k$ and $\L k' = [\L k, \L k]$ its semisimple part. This induces global decompositions $K = Z(K) \times K'$ and $K_0 = Z(K_0)K_0'$ (almost direct). Here $K'$ is a finite covering of $K_0'$, hence compact, while $Z(K)\cong \R$. In particular, both $K'$ and $Z(K)$ and hence $K$ are amenable.
\item There exists a Cartan subgroup
$H$ of $G$ with $Z(K)  \subset H \subset K$. We fix such a Cartan subgroup once and for all.
\item Denote by $\L p$ the orthogonal complement of $\L k$ in $\L p$ with respect to the Killing form so that $\L g = \L{z(k)} \+ \L k' \+ \L p$. Identify $\L p$ with the tangent space of 
symmetric space $\X$ of $G$ at the basepoint $eK$. There are two choices for the invariant complex structure of $\X$, and we fix one of them. After this choice, there exists a unique $J \in
\L{z(k)}$ such that
${\rm ad}(J)|_{\L p}$ defines the restriction of the chosen complex structure to $\L p$. 
\item Denote by $\tri = \tri(\L g_\C,
\L h_\C)$ the roots of $\L g_\C$ with respect to $\L h_\C$. Choose a positive system $\tri^+ \subset \tri$
in such a way that for all $\alpha$ in the set $\tri^+_n$ of non-compact positive roots the relation $\alpha(iJ) = 1$ holds. Fix a maximal system of strongly orthogonal roots
$\tri^{++}_n \subset \tri^+_n$.
\item The \emph{compact Weyl group} (associated to our choice of compact Cartan $H$) is defined by $W_c := N_K(H)/Z_K(H)$. This acts on $H$ by conjugation and thus on $\L h$ via the adjoint action and $\L h^*$ via the
coadjoint action.We denote by $(\L h^*)^{W_c}$ and $(\L h)^{W_c}$ the sets of $W_c$-invariants in $\L h^*$ and $\L h$ respectively.Our choice of non-compact root ensures that $\tri^+_n$ is invariant under $W_c$. 
\item Given $\alpha \in \tri$ choose root vectors $E_{\pm \alpha} \in \L g_\C^\alpha$ such that
\[i(E_\alpha + E_{-\alpha}),\; E_\alpha - E_{-\alpha} \in \L k + i \L p, \quad \alpha([E_\alpha, E_{-\alpha}]) = 2.\]
Define $h_\alpha := -i[E_\alpha,E_{-\alpha}] \in \L h$,  $X_\alpha := E_\alpha +E_{-\alpha}$ and $Y_\alpha := -i(E_\alpha-E_{-\alpha})$. Denote by $\L a$ the span of the
$X_\alpha$ for $\alpha \in \tri^{++}_n$ (which is a maximal abelian subalgebra of $\L p$) and by $A$ the associated analytic subgroup of $G$.
\item We use the isomorphism $\ker(\pi) = Z(G)  \cong \pi_1(G_0)$ to identify $\pi_1(G_0)$ with a subgroup of $G$. We observe that $\pi_1(G_0)  \cong \pi_1(K_0)$ is actually a subgroup of $K$.
\end{itemize}

\subsection{Continuous homogeneous quasimorphisms}
We keep the notations introduced in the last subsection; in particular, $G$ is a $1$-connected simple Hermitian Lie group. We describe homogeneous quasimorphisms on $G$. For background on continuous bounded cohomology and quasimorphisms we refer the reader to \cite{BC}, \cite{Bavard}, \cite{BargeGhys}. Homogeneous quasimorphisms on Hermitian Lie groups are discussed in detail in \cite{Surface}. In particular we deduce from \cite[Prop. 7.8]{Surface}:
\begin{lemma} There exists a unique homogeneous quasimorphism $ \mu_G: G  \to  \R$ satisfying $\mu_G(\exp(J)) = 1$. Any homogeneous quasimorphism on $G$ is a multiple of $\mu_G$.
\end{lemma}
The following fact could also be deduced from \cite[Prop. 7.8]{Surface}, but we prefer to give a short self-contained proof.
\begin{proposition}\label{PNeglectible} Let $\mu_G$ be a homogeneous quasimorphism on $G$. Then for all $p \in \exp(\L p)$ we have $\mu_G(p) = 0$.
\end{proposition}
\begin{proof} Since $A$ is amenable, the restriction $\mu_G|_A$ is  a homomorphism. Since homogeneous quasimorphisms are conjugation-invariant, its differential is invariant under $X \mapsto \exp({\rm ad}(\pi J))(X) = -X$. This shows that $\mu_G|_A$  is trivial, and the proposition follows by using conjugation-invariance once more.\end{proof}
By conjugation invariance the restriction $\mu_G|_K$ is uniquely determined by $\mu_G|_H$; the latter is a homomorphism, which can be determined explicitly. For this the key observation is the following lemma:
\begin{lemma}
With notation as above we have $\dim (\L h^*)^{W_c} = \dim (\L h)^{W_c} = 1$.
\end{lemma}
\begin{proof} Decompose $\L h$ into irreducibles $W_c$-modules. As a first step let $\L h' := \L h
\cap \L k' $ so that $\L h = \L{z(k)} \oplus \L h'$. Then $\L h'$ is a maximal torus in the compact semisimple Lie algebra $\L k'$, i.e.$\L h'_\C$
is a maximal torus in the complex semisimple Lie algebra $\L k'_\C$ and $W_c$ is the Weyl group associated to the pair $(\L k'_\C, \L h'_\C)$. In particular, the action of $W_c$ on
$\L{z(k)}$ is trivial, while $\L h'$ decomposes into irreducible modules corresponding to the simple subalgebras of $\L k'$. Each of these modules has dimension $\geq 3$ and is
thus non-trivial. Thus, $\dim (\L h)^{W_c} = 1$. Now fix a non-degenerate invariant bilinear form on $\L k$. The restriction of this form to $\L h$ can then be used to identify
$\L h$ and $\L h^*$ as $W_c$-modules. This yields $\dim (\L h^*)^{W_c} = \dim (\L h)^{W_c}$.
\end{proof}
Now we deduce easily:
\begin{proposition}\label{Differential}
Let notation be as above. Then for all $X \in \L h$ we have
\begin{eqnarray*}\label{FormulaDifferential} \mu_G(\exp(X)) = \frac{1}{|\tri_n^+|} \sum_{\alpha \in \tri_n^+} \alpha(iX).\end{eqnarray*}
\end{proposition}
\begin{proof} Both $d(\mu_G|_H)$ and $\sum  \alpha(i \cdot)$ define elements in $(\L h^*)^{W_c}$, hence are proportional. The proportionality constant can be computed by evaluating at $J$.
\end{proof}

\subsection{Continuous partial orders}
In this section we describe continuous partial orders on a  $1$-connected simple real Lie group $G$. We keep the notation of the last two sections.  Associated with any such order $\leq$ is a closed, topologically locally generated order semigroup $G^+$. By results of Neeb \cite{Neeb}, $G^+$ is a Lie semigroup. This means that the Lie wedge
\[C^+ := {\bf L}(G^+) = \{X\in \L g\,| \,\forall t>0: \, \exp(tX) \in G^+ \}\]
generates $G^+$ infinitesimally, i.e.
\[G^+ = \overline{\l\exp{C^+}\r}.\]
In particular, $\leq$ is uniquely determined by the ${\rm Ad}$-invariant closed, pointed
generating cone $C^+ \subset \L g$. The set $\mathcal C(\L g)$ of such cones have been described in  \cite{Vinberg} and \cite{Paneitz}. Any $C^+ \in \mathcal C(\L g)$ is determined uniquely by its intersection with $\L h$ and contains either $J$ or $-J$, and the corresponding partial order will be called \emph{positive} or \emph{negative} accordingly. Since we may always reverse the roles of $J$ and $-J$ it suffices to deal with the subset $\mathcal C(\L g)^+ \subset  \mathcal C(\L g)$ of positive orders. It turns out that not every element of $\mathcal C(\L g)^+$ is \emph{global} in the sense that it arises as the Lie wedge of an invariant continuous partial order on $G$. The subset $\mathcal C(G)^+$ of global cones in $\mathcal C(\L g)^+$ has been determined in \cite{Olshanski}. We will not need the precise classification statement, but only the following observation: Since $\mathcal C(\L g)^+ \subset  \mathcal C(\L g)$, it follows from 
\cite[Thm. 2]{Vinberg} that every $C^+ \in \mathcal C(G)^+$ contains the cone denoted $C_{\min}(\L g)$ in \cite{Vinberg}. Since the interior of $C_{\min}(\L g) \cap \L h$ contains the ray $\{t \cdot J\,|\, t > 0\}$ we deduce:
\begin{lemma}\label{TT}
Let $C^+ \in \mathcal C(G)^+$ and put $c := C^+ \cap \L h$. Then the interior $c^\circ $ of $c$ in $\L h$ contains $\{t \cdot J\,|\, t > 0\}$.
\end{lemma} 

\section{Realizing quasimorphisms as relative growth}\label{SecMain}

\subsection{The sandwich condition}

We have observed in the introduction that a quasimorphism can be reconstructed as relative growth from the associated partial order. A more general condition allowing for such a  reconstruction is the following:

\begin{definition}\label{DefAssoc} A non-zero homogeneous quasimorphism $f$ on a topological group $G$ is said to \emph{sandwich} an invariant partial order $\leq$ if there exist constants $C_1, C_2 \in \R$ such that
\begin{eqnarray}\label{Sandwich0}Q^+_f(C_1) \subset G^+ \subset Q^+_f(C_2),\end{eqnarray}
where for $C \in \R$ the  \emph{superlevel set} $Q^+_f(C) $ is given by.
\[Q^+_f(C) := \{g \in G\,|\, f(g) \geq C\}\]
\end{definition}
In fact, the upper bound comes for free:
\begin{lemma}\label{NoUpperBound}Let $G$ be group, $f: G \to \R$ a non-trivial homogeneous quasimorphism and $\leq$ be an invariant partial oder on $G$ with order semigroup $G^+$. If there exists $C_1>0$ with 
\[Q^+_f(C_1) \subset G^{+},\]
then $\leq$ is sandwiched by $f$, $Q^+_f(C_1) \subset G^{++}$, and \eqref{Sandwich0} is satisfied with $C_2 := 0$.
\end{lemma}
\begin{proof} First we claim that $G^+ \subset Q^+_f(0)$. Otherwise we find $g_0 \in G^+$ with $f(g_0) < 0$. Then every $g \in G$ can be written as $g = g_0^{n}(g_0^{-n}g)$, where $n \in \mathbb N$ is chosen in such a way that $f(g_0^{-n}g) > C$. This implies $g_0^{-n}g \in G^+$ and thus $g \in G^+$. Since $g \in G$ was arbitrary this implies $G^+ = G$ contradicting the pointedness of $G^+$. This proves our claim. Now assume $Q^+_f(C) \subset G^{+}$ and suppose $g \in G$ satisfies $f(g) \geq C >0$. Then for any $h \in G$ we find $n \in \mathbb N$ such that $f(g^nh^{-1}) > C$, whence $g^n > h$, showing that already $g \in G^{++}$.
\end{proof}
The following generalization of Proposition \ref{dull} will be at the heart of our proof of Theorem \ref{IntroMain}; we therefore give its elementary proof in some details.
\begin{proposition}\label{Sandwich}
Suppose that $(G, \leq)$ is ordered group and that $f: G \ra \R$ is a non-zero homogeneous quasi-morphism. If $f$ sandwiches $\leq$, then $\leq$ is admissible and
for all $g \in G^{++}$, $h \in G$ we have
\[\gamma(g,h) = \frac{f(h)}{f(g)}.\]
\end{proposition}
\begin{proof} Let $g \in G^{++}$ and $h \in G$. Define $T_n(g,h):= \{p \in \Z\,|\, g^p \geq h^n\}$ and $\gamma_n(g,h) = \inf T_n(g,h)$ so that
\begin{eqnarray}\label{GammaN}
\gamma(g,h) = \lim_{n \to \infty}\frac{\gamma_n(g,h)}{n}.
\end{eqnarray}
Choose a constant $C_1 >0 $ such that $Q^+_f(C_1) \subset G^{++}$ holds. Since $f$ is non-trivial, $Q^+_f(C_1)$ is non-empty, and thus $\leq$ is admissible. We claim that any
integer $p_n$ satisfying
\[p_n \geq \frac{n f(h) + C_1 + D(f)}{f(g)}\]
also satisfies  $p_n \in T_n(g,h)$. (Such a $p_n$ exists since $f(g) \neq 0$ for any dominant $g$.) Indeed, we have
\[f(g^{p_n}h^{-n}) \geq p_n f(g) - nf(h) - D(f) \geq (n f(h) + C_1 + D(f)) - nf(h) - D(f) = C_1,\]
hence $g^{p_n}h^{-n} \in G^+$, which implies $g^{p_n} \geq h^n$ as claimed. In particular $\gamma_n(g,h) \leq p_n$ and choosing $p_n$ minimal possible we obtain
\begin{eqnarray}\label{MainUpperBound}
\frac{\gamma_n(g,h)}{n} \leq \frac{f(h)}{f(g)} + \frac{C_1 + D(f) + f(g)}{nf(g)}.
\end{eqnarray}
Now suppose $p \in \Z$ satisfies
\[p < \frac{nf(h) - D(f)}{f(g)}.\]
Then
\[f(g^{p}h^{-n}) \leq pf(g) -nf(h) + D(f) < (nf(h) - D(f)) -nf(h) + D(f) = 0.\]
Thus $g^{p}h^{-n} \not \in G^+$ and thus $p \not \in T_n(g,h)$. Consequently,
\begin{eqnarray}\label{MainLowerBound}
\frac{\gamma_n(g,h)}{n} \geq \frac{f(h)}{f(g)} - \frac{D(f)}{nf(g)}.
\end{eqnarray}
Combining \eqref{MainUpperBound} and \eqref{MainLowerBound} and passing to the limit $n \to \infty$ we obtain the proposition.
\end{proof}
We now turn to the proof of Theorem \ref{IntroMain}. By Proposition \ref{Sandwich} and Lemma \ref{NoUpperBound} it suffices to establish the following estimate:
\begin{lemma}\label{MainLemma} Let $\leq$ be a continuous partial order with order semigroup $G^+$ on a $1$-connected Hermitian simple Lie group $G$. Then there exists $C>0$ such that
\begin{eqnarray}\label{QuasimorphismVersusCone}  Q^+(C) := Q^+_{\mu_G}(C) \subset G^+.\end{eqnarray}
\end{lemma}
The remainder of this section is devoted to the proof of Lemma \ref{MainLemma} (and hence Theorem \ref{IntroMain}). We claim that Lemma \ref{MainLemma} can be deduced from the following two lemmata:
\begin{lemma}\label{CompactLemma}
There exists a constants $C_1 \in \R$ such that
\begin{eqnarray}\label{CompactBounds} Q_+(C_1) \cap K \subset G^+ \cap K.\end{eqnarray}
\end{lemma}
\begin{lemma}\label{NoncompactLemma}
There exists a constant $C_0$ such that for all $p \in \exp_{G}(\L p) \subset G$ there exists $k(p) \in K$ with $|\mu_G(k(p))| \leq C_0$ and $k(p) p \geq e$.
\end{lemma}
Before we prove the lemmata, let us check carefully that they imply Lemma \ref{MainLemma}. 
Every $g \in G$ can be written as $g = kp$ with $k \in K$, $p \in \exp_G(\L p)$. According to Lemma
\ref{NoncompactLemma} we can choose $k(p)$ with $k(p)p \geq e$ and  $|\mu_G(k(p))| \leq C_0$. Now we claim that Lemma \ref{CompactLemma} provides the desired estimate \eqref{QuasimorphismVersusCone} for
\[C := C_0 + C_1 +D(\mu_G).\]
Indeed, suppose $\mu(g) \geq C$. Then using Proposition \ref{PNeglectible} and \eqref{CompactBounds} we deduce 
\begin{eqnarray*}
&&\mu_G(k) \geq \mu_G(kp) - \mu_G(p) - D(\mu_G) = \mu_G(g) -D(\mu_G) \geq C_0 + C_1\\
&\Rightarrow& \mu_G(kk(p)^{-1}) = \mu_G(k) - \mu_G(k(p)) \geq C_1\\
&\Rightarrow& kk(p)^{-1} \in Q_+(C_1') \cap K \subset G^+ \cap K\\
&\Rightarrow& k \geq k(p)\\ &\Rightarrow& g = kp \geq k(p)p \geq e.
\end{eqnarray*}
This proves the claim and reduces Theorem \ref{IntroMain} to the above two lemmata, whose respective proofs will be the content of the following two subsections.

\subsection{Proof of the main theorem I: The compact contribution}

In this subsection we prove Lemma \ref{CompactLemma}. Thus assume $k \in Q_+(C_1) \cap K$ for some $C_1 > 0$. Then $k$ is conjugate to some $h \in H$ with $\mu_G(h) > C_1$. Since $G^+$ is conjugation-invariant we have $k \in G^+$ iff $h \in G^+$. We have thus reduced Lemma \ref{CompactLemma} to the following observation:
\begin{lemma}
There exists a constant $C_1 \in \R$ such that
\begin{eqnarray}\label{HBounds} Q_+(C_1) \cap H \subset G^+ \cap H.\end{eqnarray}
\end{lemma}
\begin{proof} We decompose $\L h = \L{z(k)}\oplus \L h'$; accordingly, every $h \in H$ can be written as 
\begin{eqnarray}\label{StandardFormH} h = \exp(tJ + X) = \exp(tJ)\exp(X) \quad (t \in \R, X \in \L h').\end{eqnarray}
Here $\exp(X) \in \exp(\L h')$, which is a compact group. Since $\mu_G$ restricts to a homomorphism on $H$ and there are no non-trivial homomorphisms from a compact group into $\R$ we have $\mu_G(\exp(X)) = 0$ and hence
\begin{eqnarray*}\label{StandardFormmu}\mu_G(h) = \mu_G(\exp(tJ)) = t.\end{eqnarray*}
Denote by $\Lambda$ the kernel of the exponential function $\exp: \L h' \to H$. Then $\Lambda$ is a cocompact lattice in $\L h'$ and thus has a bounded fundamental domain in $\L h'$. Consequently, if we denote by $B'_r(0)$ the closed ball of radius $r$ around $0$ in $\L h'$ then 
\begin{eqnarray}\label{Minkowski}
\exists R > 0 \;\forall r \geq R \;\forall X \in \L h' \;\exists Y \in \Lambda: X + Y \in B'_r(0).
\end{eqnarray}
Now let $c := {\bf L}(G^+) \cap \L h$ and denote by $c^\circ$ the interior of $c$. By Lemma \ref{TT} we have $tJ \in c^\circ$ for all $t>0$. Denote by $S_t := tJ + \L h'$ the affine hyperplane through $tJ$ parallel to $\L h'$. Then $S_t \cap c^\circ$ is open in $S_t$ and non-empty, since it contains $tJ$. We thus see that for $t>0$ we have
\begin{eqnarray*}\label{Radius}r(t) := \max\{r > 0\,|\, tJ + B'_{r}(0)\subset c\} > 0.\end{eqnarray*}
In fact, convexity of $c$ implies $r(t) \to \infty$ as $t \to \infty$. By \eqref{Minkowski} we thus have
\begin{eqnarray}\label{Mink2}\exists T > 0 \;\forall t \geq T \;\forall X \in \L h \;\exists Y \in \Lambda: tJ + X + Y \in c.\end{eqnarray}
Now we claim that we can choose $C_1 := T$. Indeed, if $h$ is as in \eqref{StandardFormH} and $\mu_G(h) = t \geq T$, then by \eqref{Mink2} we find $Y \in \Lambda$ with $tJ + X + Y \in c$ and hence
\[h = \exp(tJ + X) = \exp(tJ + X + Y) \in \exp(c) \subset H \cap G^+.\]
\end{proof}
This finishes the proof of Lemma \ref{CompactLemma}. 

\subsection{Proof of the main theorem II: The non-compact contribution} The purpose of this subsection is to establish Lemma \ref{NoncompactLemma}, thereby finishing the proof of Theorem \ref{IntroMain}. We will argue by reduction to the case of the universal covering group of $SL_2(\R)$, which we denote by $\widetilde{SL}_2(\R)$. (This case was treated in \cite{Gabi}.) We recall that for any $\alpha \in \tri_n^{++}$ the bracket relations $[X_\alpha, Y_\alpha] = -2
h_\alpha$, $[h_\alpha, X_\alpha] = 2Y_\alpha$ and $[h_\alpha, Y_\alpha] = -2X_\alpha$ hold. (See e.g.\cite{Paneitz}, where the notation is compatible with ours.)
Therefore, the three-dimensional real Lie algebra $\L{sl}_\alpha$ spanned by $X_\alpha, Y_\alpha, h_\alpha$ is isomorphic to $\L {sl}_2(\R)$ via an isomorphism $\sigma_\alpha$ given by
\[\sigma_\alpha(X_\alpha) := \begin{pmatrix} - 1 & 0 \\0 & 1\end{pmatrix}, \, \sigma_\alpha(Y_\alpha) :=  \begin{pmatrix} 0 & 1 \\-1 & 0\end{pmatrix}, \, \sigma_\alpha(h_\alpha) := \begin{pmatrix} 0 &
1 \\1 & 0\end{pmatrix}.\]
Denote by $\psi_\alpha: \L{sl}_2(\R) \to \L{sl}_\alpha \hookrightarrow \L g$ the inclusion induced by the inverse of this isomorphism. Then $\psi_\alpha$ integrates to a group
homomorphism
\[\Psi_\alpha:\widetilde{SL}_2(\R) \to G.\]
In fact, $\Psi_\alpha$ factors through a map $\Psi_\alpha^0:SL_2(\R) \to G_0$, in particular \begin{eqnarray}\label{FundamentalGroups}\Psi_\alpha(\pi_1(SL_2(\R))) \subset \pi_1(G_0).\end{eqnarray} Indeed, since $SL_2(\C)$ is simply-connected the complexification \[(\psi_\alpha)_\C: \L{sl}_2(\C) \to \L g_\C\] integrates to a map $(\Psi_\alpha)_\C: SL_2(\C) \to (G_0)_\C$. Since $G_0$ is
linear and connected, it coincides with the analytic subgroup of its universal complexification with Lie algebra $\L g$ \cite[Satz I.6.1 and Satz III.9.24]{HiNe}. Now $(\psi_\alpha)_\C$ maps
$\L{sl}_2(\R)$ into $\L g$ and thus $(\Psi_\alpha)_\C$ maps $SL_2(\R)$ into $G_0$. Then the restriction of $(\Psi_\alpha)_\C$ provides the desired factorization $\Psi_\alpha^0$ of $\Psi_\alpha$. We will now use the maps $\Psi_\alpha$ to reduce our problem to the case of $\widetilde{SL}_2(\R)$ by means of the following lemma:
\begin{lemma}\label{PositivityPreserved}
For each $\alpha \in \tri_n^{++}$ there exists a continuous admissible partial ordering $\leq_\alpha$ on $\widetilde{SL}_2(\R)$ with the following property: If
$g \geq_{\alpha} e$ for some $g \in \widetilde{SL}_2(\R)$, then $\Psi_\alpha(g) \in G^+$.
\end{lemma}
\begin{proof} Denote by $C^+$ the Lie wedge of $G^+$ and define \[C_\alpha^+ := \sigma_\alpha(C^+ \cap \L{sl}_\alpha) \subset \L{sl}_2(\R).\]
Since the kernel of the map
\[\Psi_\alpha: \widetilde{SL}_2(\R) \to \Psi_\alpha(\widetilde{SL}_2(\R))\]
is central, $\Psi_\alpha$ induces an isomorphism $\rm{Ad}(\widetilde{SL}_2(\R)) \to {\rm Ad}(\Psi_\alpha(\widetilde{SL}_2(\R)))$. As $\sigma_\alpha$ is equivariant with respect
to these adjoint actions, we deduce that $C^+_\alpha$ is an $\rm{Ad}$-invariant, closed pointed cone in $\L{sl}_2(\R)$. This cone is non-trivial, since $h_\alpha \in C^+$, and thus $\L{sl}_2(\R) = C^+_\alpha - C^+_\alpha$, since  the right hand side is a non-trivial ideal. This means that $C^+_\alpha$ is generating. Now there exists only two (mutually inverse) $\rm{Ad}$-invariant, closed pointed generating cones in $\L{sl}_2(\R)$, and both are global. This means that there exists a partial order $\leq_\alpha$ on $\widetilde{SL}_2(\R)$ with order semigroup
$\overline{\langle \exp(C^+_\alpha)\rangle}$. Since \[\psi_\alpha(C^+_\alpha) = C^+ \cap \L{sl}_\alpha \subset C^+\] we have $\Psi_\alpha(\exp(C^+_\alpha)) \subset G^+$, from
which the lemma follows.
\end{proof}
To finish our argument we use the following fact about the $\widetilde{SL}_2(\R)$-case:
\begin{lemma}\label{ReduceToMe}
Given any continuous admissible ordering $\leq$ on $G = \widetilde{{\rm SL}_2}(\R)$ there exists a constant $N$ and an element $z_0 \in \pi_1({\rm SL}_2(\R))$ such that
for every $X \in \L p$ there exists $0 \leq n \leq N$ with
\[z_0^n \exp_G(X) \geq e.\]
\end{lemma}
\begin{proof} In \cite{Gabi} a continuous admissible positive order on $G$ was introduced by dynamical means. By \cite[Sec. 3.5]{Vinberg} this is the only continuous admissible ordering on $G$ up to inversion, and we denote it by $\leq$. We know from \cite[Lemma 2.9]{Gabi} (specialized to $n=1$) that for every $X \in \L p$
there exists a path $g_t$ defining a homotopy class $g = [g_t] \in G$ which satisfies both $g_1 = \exp_{SL_2(\R)}(X)$ and $g \geq e$ and has Maslov quasimorphism $\mu_{Maslov}(g) \leq 4 \pi$. If we denote $p := \exp_G(X)$, then the first statement means that $g = zp$ for some $z \in \pi_1({\rm SL}_2(\R)) \subset G$. If $z_0$ denotes the positive generator of $\pi_1({\rm SL}_2(\R))\cong \Z$, then $z = z_0^n$ for some $n > 0$, and the uniform bound on the Maslov quasimorphism implies the uniform bound on $n$. This establishes the lemma for $\leq$, hence for all continuous admissible orderings on $G$.
\end{proof}
Now we can finish the proof of Lemma \ref{NoncompactLemma}: Combining Lemma \ref{PositivityPreserved} and Lemma \ref{ReduceToMe} we now choose for every $\alpha \in \tri_n^{++}$ a constant $N_\alpha \in \mathbb N$ and an element $z_{\alpha,0} \in \pi_1({\rm SL}_2(\R))$ such
that for every $t_\alpha \in \R$ there exists $0 \leq n_\alpha \leq N_\alpha$ with
\[z_{\alpha,0}^{n_\alpha} \exp\left(t_\alpha \begin{pmatrix} - 1 & 0 \\0 & 1\end{pmatrix}\right) \geq_\alpha e.\]
Define $z_\alpha := \Psi_\alpha(z_{\alpha,0})$. By \eqref{FundamentalGroups} we have $z_\alpha \in \pi_1(G_0) = Z(G)$. Applying $\Psi_\alpha$ and using Lemma
\ref{PositivityPreserved} we obtain:
\[z_\alpha^{n_\alpha} \exp(t_\alpha X_\alpha) \in G^+.\]
Now, any $a \in A$ is of the form
\[a = \prod_{\alpha \in \tri^{++}} \exp(t_\alpha X_\alpha)\]
for some $t_\alpha \in \R$, and any $g \in \exp(\L p)$ is of the form $g=kak^{-1}$ for some $k \in K$. This implies that for every $g \in \exp(\L p)$ we can find $0 \leq n_\alpha \leq
N_\alpha$ such that
\[\left(\prod_{\alpha \in \tri_n^{++}} z_\alpha^{n_\alpha}\right) \cdot g = k \left(\prod_{\alpha \in \tri^{++}} z_\alpha^{n_\alpha} \exp(t_\alpha X_\alpha)\right)  k^{-1} \in
G^+.\]
This implies Lemma \ref{NoncompactLemma} and finishes the proof of Theorem \ref{IntroMain}.
\subsection{Beyond simple-connectedness}\label{NonSC}
In the proof of Theorem \ref{IntroMain} we have always assumed $G$ to be simply-connected. This assumption ensured in particular the existence of a non-zero homogeneous quasimorphism and a
non-trivial continuous admissible partial order on $G$. As far as the former existence question is concerned, it is easy to classify the non-simply connected simple Lie groups
$\widehat{G}$ which
admit a non-zero homogeneous quasimorphisms. For this we recall that the space of such quasimorphisms is
\[EH^2_{cb}(\widehat{G}; \R) = \ker(H^2_{cb}(\widehat{G}; \R) \to H^2_{c}(\widehat{G}; \R)).\]
Since $\dim H^2_{cb}(\widehat{G}; \R) \leq 1$ this is equivalent to $H^2_{cb}(\widehat{G}; \R) \cong \R$ and $H^2_{c}(\widehat{G}; \R) = 0$. Equivalently, $\widehat{G}$ is
Hermitian with finite fundamental group. Denote by $p: G \to \widehat{G}$ its universal covering. Then there is a unique homogeneous quasimorphism $\mu_{\widehat{G}}: \widehat{G} \to \R$ such that $p^*\mu_{\widehat{G}} = \mu_G$. Now it follows from the fact that $\ker(p)$ is torsion that for every continuous order on $G$ with order semigroup $G^+$ the image $p(G^+)$ is again the order semigroup of a continuous order on $\widehat{G}$, and these are actually all continuous orders on $\widehat{G}$. It then follows immediately that every continuous order on $\widehat{G}$ is sandwiched by either $\mu_{\widehat{G}}$ or $-\mu_{\widehat{G}}$. This implies that Theorem \ref{IntroMain} holds in fact for all simple Hermitian Lie groups with finite fundamental group.

\section{Implications and further examples}\label{SecApps}

\subsection{Basic definitions}
Among the initial motivation of Eliashberg and Polterovich to introduce relative growth was the construction of a certain metric $G$-space out of an admissible ordered group $(G,\leq)$. To explain their construction, let $G$ be a group and $\leq$ an admissible invariant ordering on $G$. Then the restriction of the relative growth function defines a positive function
\[\gamma: G^{++} \times G^{++} \to \R^{>0},\]
whose symmetrized logarithm
\[d(g,h) := \log \max\{\gamma(g,h), \gamma(h,g)\}\]
yields a pseudo-metric on $G^{++}$. We refer to the associated metric space as the \emph{order space} of $(G, \leq)$ and denote it by $\mathfrak X(G, \leq)$. Note that the conjugation action of $G$ on $G^{++}$ induces an isometric $G$-action on $\mathfrak X(G, \leq)$. In general, it is a difficult problem to compute the order space of an ordered group. However, if the order in question is sanwiched by a homogeneous quasimorphism, then we can apply Proposition \ref{Sandwich} in order to compute the order space explicitly:
\begin{corollary}\label{GeneralCollapse}
Suppose that $(G, \leq)$ is an admissible ordered group and that $f: G \ra \R$ is a continuous homogeneous quasi-morphism sandwiching $\leq$. Then the map
\[\iota: \mathfrak X(G,\leq) \ra \R, \quad [g] \mapsto \log f(g)\]
is an isometry onto its image.
\end{corollary}
\begin{proof} Let $g,h \in G^{++}$. By Proposition \ref{Sandwich} we have
\[d([g],[h]) = \max \{\log \gamma(g,h), \log \gamma(h,g)\} = |\log f(g) -\log f(h)|,\]
showing that $\iota$ is an isometry.
\end{proof}

\subsection{Finite-dimensional examples}
Applying Corollary \ref{GeneralCollapse} to the case of $1$-connected, finite-dimensional simple Lie groups discussed in the main theorem we obtain:
\begin{corollary}\label{MainCorollary}
\begin{itemize}
\item[(i)] Let $G$ be a 1-connected simple Hermitian Lie group equipped with an arbitrary continuous order $\leq$. Then there is a surjective isometry
\[\iota: \mathfrak X(G,\leq) \ra \R, \quad [g] \mapsto \log \mu_G(g).\]
\item[(ii)] Denote by $\leq$ the admissible ordering on $K$ obtained by restricting a continuous order from $G$. Then there is still a surjective isometry
\[\iota: \mathfrak X(K,\leq) \ra \R, \quad [k] \mapsto \log \mu_G(k).\]
\end{itemize}
\end{corollary}
Indeed, Corollary \ref{GeneralCollapse} applies in view of Lemma \ref{MainLemma}  and Lemma \ref{CompactLemma} respectively, and surjectivity follows from $\mu_G(\exp(tJ)) = t$ in both cases. In fact it is easy to see that every order space of a Lie group necessarily contains a copy of $\R$ as the image of a suitable one-parameter semigroup in $G^{++}$. In that sense the order spaces of 1-connected simple Hermitian Lie group with respect to a given continuous order is as small as possible for a Lie group. Corollary \ref{MainCorollary} answers a question of Polterovich, which was the starting point for the investigations in this paper. The results in Corollary \ref{MainCorollary} should be compared to the case of $1$-connected, finite-dimensional \emph{abelian} Lie groups, i.e. finite-dimensional vector spaces.
\begin{example}\label{LinearOrder}
Let $V$ be a finite-dimensional vector space (considered as an abelian Lie group under addition) and $C^+ \subset V$ a closed, pointed convex cone with non-empty interior. 
By \cite[Corollary 11.7.1]{Convexity} there exists a weak-$*$-compact subset $\mathcal A^*$ of the unit ball $V^*_1$ in $V^*$ such that
\begin{eqnarray}\label{hyperplanes} C^+ = \{v \in V\,|\,\forall\alpha \in \mathcal A^*: \, \alpha(v) \geq 0\}.\end{eqnarray}
The dominants of the partial order with order semigroup $C^+$ are given by $C^{++} = {\rm Int}(C^+)$.  A short computation shows that the pseudo-distance $d$ on $C^{++}$ is given by \begin{eqnarray}\label{LinearDistance}d(v,w) = \max_{\alpha \in \mathcal A^*} |\log\alpha(v) - \log \alpha(w)|.\end{eqnarray} This is actually a metric on $C^{++}$, and thus $\mathfrak X(V, \leq) = ( {\rm Int}(C^+), d)$.
\end{example}
Thus in the abelian case, the order space is as large as possible (i.e. the natural map $G^{++} \to \mathfrak X(G, \leq)$ is one-to-one), while in the simple case it is as small as possible.

\subsection{Infinite-dimensional examples}

The strong dichotomy between order spaces of finite-dimensional simple and finite-dimensional abelian Lie groups discovered in the last subsection exists also for certain families of infinite-dimensional Lie groups, which we discuss here. For this we return to the original setup, in which relative growth was introduced, namely contact and symplectic geometry. Various infinite-dimensional Lie groups with natural invariant orders arise in this context, and for several classes of such groups the associated order spaces have been studied in  \cite{EP, Gabi}. In all these examples the order spaces turn out to be infinite-dimensional. In this subsection we provide an example of a similar geometric flavour, in which the order space fails not only to be infinite-dimensional, but in fact collapses to $\R$. The reason for this collapse is again provided by a homogeneous quasimorphism, which sandwiches the order in question.\\

In order to explain our example, we introduce the following notation: Denote by $(M,\omega)$ a closed symplectic manifold of dimension $2n$. Every
smooth, time dependent function $H_{t}:M\rightarrow \mathbb{R}$ gives rise to a smooth vector field, $X_{H}$ via the
pointwise linear equation $dH=-\omega(X_{H},\cdot)$. These vector fields are called \emph{Hamiltonian}. The group $G_0 := {\rm Ham}(M, \omega)$ of \textit{Hamiltonian motions} is by definition the subgroup of the diffeomorphism
group ${\rm Diff}(M)$ given by the time-1 maps of the flows generated by the Hamiltonian vector fields. Since $\omega^n$ is a volume form on $M$, $G_0$ is actually a subgroup of the volume preserving
diffeomorphisms of $M$. A detailed study of the group $G_0$ is provided in \cite{Polt}. Here we just remark that $G_0$ admits a natural topology and smooth structure, turning it into an
infinite-dimensional Lie group. We will be interested in the universal covering $G$ of $G_0$.\\

An important problem is the existence and uniqueness problem for Calabi type quasimorphisms on $G$. For background on this complex of problems see
\cite[Chapter 10]{McDuff}. Here we recall only some of the most basic definitions in order to fix our notation: Given an open subset $U \subset M$, denote by $G_U^0$ the group of
Hamiltonian diffeomorphisms of $M$ generated by Hamiltonians supported inside $U$, and observe that the elements of $G_U^0$ are then automatically compactly supported. On the
universal covering $G_U$ of $G_U^0$ there exists a homomorphism ${\rm Cal}_U: G_U \to \R$ called the \emph{Calabi homomorphism} given as follows: If $[f_t] \in G_U$ is
represented by a path $f_t$ in $G_U^0$ generated by a time-dependent Hamiltonian $F_t$, then
\[{\rm Cal}_U([f_t]) = \int_0^1 \int_U F_t\omega^{n}dt.\]
This homomorphism descends to a homomorphism of $G_U^0$ if $\omega|_U$ is exact, but not in general. Now let us call an open subset $U \subset M$ \emph{displacable} if there
exists $g \in G_0$ such that $gU \cap \overline{U} = \emptyset$. Then we define:

\begin{definition}
A quasimorphism $f: G \to \R$ is called \emph{of Calabi type} if for every displacable open subset $U \subset M$ the equality
\[f|_U = {\rm Cal}_U\]
holds.
\end{definition}

In \cite{EntP} the existence of a Calabi type quasimorphism was established for symplectic manifolds which are spherically monotone and the even part of whose quantum homology
algebra is semisimple. We cannot explain these assumptions here, but refer the reader to the aforementioned article and the references therein for details. Here we can only sketch some
ideas of the construction. The basic idea of Entov and Polterovich for constructing a Calabi type quasimorphism is to use the \emph{spectral invariants} of $G$. For the present purpose
it suffices to know that these are given by a map
\[c: QH_{ev}(M) \times G \to \R, \quad (a,g) \mapsto c(a,g),\]
where $QH_{ev}(M)$ is the even part of the quantum homology algebra of $M$. Then they prove the following result:

\begin{lemma}[Entov-Polterovich]\label{Ent}
If $QH_{ev}(M) = Q_1 \oplus \dots \oplus Q_d$ denotes the decomposition of $QH_{ev}(M)$ into a direct sum of fields and $e$ the unit of $Q_1$, then \[r := - c(e, \cdot): G \to \R\] is a continuous quasimorphism. Moreover,
\begin{itemize}
\item[(i)] $r(gh) \geq r(g) + r(h)$  for all $g,h \in {G}$;
\item[(ii)] If $e_G$ denote the identity element of $G$, then $r(e_G) = 0$;
\item[(iii)] $r$ is conjugation-invariant.
\end{itemize}
Up to a constant factor of ${\rm Vol}(M)$ the homogeneization $\widetilde{r}$ of $r$ is of Calabi type.
\end{lemma}

For proofs see again \cite{EntP}, in particular Section 2.6, Theorem 3.1 and Proposition 3.3. We refer to $r$ as the \emph{spectral quasimorphism} on $G$. Based on the examples
from finite-dimensional Lie groups it is reasonable to ask whether $\widetilde{r}$ sandwiches a partial order. For the study of this problem, we suggest the following terminology:

\begin{definition}
A closed symplectic manifold $(M,\omega)$ is called \emph{Calabi orderable} if there exists a partial order $\leq$ on $G$ and a dominant $g \in G^{++}$ with respect to this
ordering such that the relative growth $\gamma(g, \cdot)$ is a Calabi type quasimorphism. In this case, $\leq$ is called a \emph{Calabi order} on $G$.
\end{definition}

We will now provide criteria which guarantee Calabi orderability. We call the spectral quasimorphism $r$ \emph{non-degenerate} if it satisfies
\[r(g) = r(g^{-1}) = 0 \Rightarrow g = e_G\]
for all $g \in G$. In this situation, Lemma \ref{Ent} yields immediately the following corollary:

\begin{corollary}\label{SpectralOrder}
Let $(M,\omega)$ be a spherically monotone closed symplectic manifold the even part of whose quantum homology
algebra is semisimple and whose spectral quasimorphism is non-degenerate. Then
\begin{itemize}
\item[(i)]  The set
\[G^+ := \{g \in G\,|\, r(g^{-1}) \leq 0\}\]
is a closed, conjugation invariant pointed submonoid of $G$ and thus defines a partial order $\leq$ on $G$.
\item[(ii)] The homogeneization $\widetilde{r}$ (and hence the Calabi quasimorphism $\widetilde{\mu} := {\rm Vol}(M) \cdot \widetilde{r}$) sandwich $\leq$.
\end{itemize}
\end{corollary}
We refer to the order $\leq$ from Proposition \ref{SpectralOrder} as the \emph{spectral order} on $G$. We briefly recall some conditions that guarantee the non-degeneracy of the spectral quasimorphism:
\begin{definition} A closed symplectic manifold $(M, \omega)$ is called
\begin{itemize}
\item \emph{rational} if $\omega(\pi_2(M)) \subset \R$ is a discrete subset;
\item \emph{strongly semipositive}, if there is no spherical homology class $A \in \pi_2(M)$ such that $\omega(A) > 0$ and $2 - n \leq c_1(A) < 0$.
\end{itemize}
\end{definition}
Then we have:
\begin{theorem} Let $(M,\omega)$ be a spherically monotone closed symplectic manifold the even part of whose quantum homology
algebra is semisimple. If $M$ is rational and strongly semipositive, then it is Calabi orderable. More precisely, a Calabi order is given by the spectral order $\leq$. Moreover,
$\mathfrak X(G, \leq) \cong \R$.
\end{theorem}
\begin{proof} By a result of Oh \cite[Theorem A]{Oh} the conditions on $M$ ensure that the spectral quasimorphism is non-degenerate. Thus the Calabi type quasimorphism
$\widetilde{\mu}$ of Entov and Polterovich sandwiches the spectral order and the result follows. 
\end{proof}
The theorem applies in particular to $\C\mathbb P^n$ with the Fubini-Study form; in particular
\[\mathfrak X(\widetilde{\rm Ham}(\C\mathbb P^n), \leq) \cong \R\]
is \emph{not} infinite-dimensional.

\subsection{Collapse of the order space in the absence of quasimorphisms}
We have seen examples of both finite- and infinite-dimensional ordered Lie groups for which the order space is much smaller than expected. This collapsing phenomenon could in both cases be tracked back to the existence of a certain homogeneous quasimorphisms and one might thus get the impression that homogeneous quasimorphisms are the only reason for a collapse of the order space. The following example shows that this is not the case and, in fact, that the order space can collapse even in the complete absence of quasimorphisms: Consider the standard embeddings
\[{\rm Sp}(2, \R) \subset {\rm Sp}(4, \R) \subset {\rm Sp}(6, \R) \subset \dots\]
induced from the embeddings
\[T^*\R \subset T^*\R^2 \subset  T^*\R^3 \subset \dots\]
Let us abbreviate by $G_n$ the universal covering of ${\rm Sp}(2n, \R)$. Then we have a similar chain for the groups $G_n$. Each $G_n$ carries a unique continuous admissible partial ordering (the maximal partial ordering in the notation of the last section) and we denote the associated order semigroup by $G_n^+$. We then define
\[G := \lim_{\rightarrow} G_n = \bigcup G_n, \quad G^+ := \bigcup G_n^+ \subset G.\]
It turns out that $G^+$ defines an admissible order $\leq$ on $G$ with
\[G^{++} =\bigcup G_n^{++} .\]
We claim that \begin{eqnarray}\label{CollapseLimit} \mathfrak X(G,\leq) \cong \R.\end{eqnarray} Indeed, let $J \in C^{++}_1$ be the element in the center of $\L k$ defining the complex strucure on the symmetric space and $L := \{\exp(tJ)\,|\, t >0\} \subset G^{++}$. We
denote by $[L]$ the image of $L$ in $\mathfrak X(G, \leq)$. Clearly, $[L] \cong \R$. We
claim that $\mathfrak X(G, \leq) = [L]$. Indeed, let $g \in G^{++}$ and choose $n \in \mathbb N$ such that $g \in G_n$. Since $\exp(J)$ and $g$ are dominant in $G_n$ we have both $\mu_{G_n}(\exp(J)) >
0$ and $\mu_{G_n}(g) >0$. Thus, there exists $t > 0$ such that $\mu_{G_n}(\exp(tJ)) = \mu_{G_n}(g)$. We now consider $\exp(tJ)$ and $g$ as elements of $G_n^{++}$ and
denote by $d_{G_n}(\exp(tJ), g)$ the corresponding pseudo-distance. By Corollary \ref{MainCorollary} we have
\[d_{G_n}(\exp(tJ), g) = 0.\] This implies \[d_G(\exp(tJ), g) = 0, \] since the natural map $\mathfrak X(G_n, \leq) \to \mathfrak X(G, \leq)$ is contractive.
This shows that $[g] = [\exp(tJ)] \in \mathfrak X(G, \leq)$. Since $\exp(tJ) \in L$ we have $[g] \in [L]$ as claimed. This establishes \eqref{CollapseLimit}. On the other hand, applying the Kotschick swindle \cite{Kotschick} to the diagonal $\widetilde{SL}_2(\R)$-subgroups of $G$, we see immediately that every homogeneous quasimorphism $f$ on $G$ restricts to a homomorphism on $G_1$. Since $G_1$ is simple, this homomorphism is trivial. But the only quasimorphism of $G_n$ restricting trivially to $G_1$ is the trivial one, whence $f$ must be trivial on every $G_n$, hence on $G$. This shows that $G$ does not possess any non-trivial homogeneous quasimorphism.

\subsection{Open problems}
We have seen in various examples that suitable homogeneous quasimorphisms allow the explicit computation of relative growth and, consequently, order spaces for ordered Lie groups. There are various directions into which our results can be extended. As far as finite-dimensional Lie groups are concerned, we have dealt with the extremal cases of simple and abelian Lie groups. In view of the structure theory of finite-dimensional Lie groups, the next step towards a complete understanding of order spaces would be to understand the behaviour of relative growth under semidirect products. For non-semisimple Lie groups the order space will probably not collapse, since the quasimorphism becomes trivial on the radical, whence it would be interesting to compute its precise form.\\

A second direction to be pursued is obviously the study of infinite-dimensional Lie groups. Here the interest is probably not in maximal generality, but rather in concrete computations of relative growth for specific classes of groups arising in contact and symplectic geometry. In the finite-dimensional case, a key step towards our computations was the reduction of continuous orders to invariant cones inside the Lie algebra. It would be interesting to know, whether such a reduction can also be used in the infinite-dimensional context.

\def\cprime{$'$}

\bigskip

\textsc{Departement Mathematik, ETH Z\"urich, R\"amistrasse 101, 8092 Z\"urich, Switzerland}\\
\texttt{$\{$gabi.ben.simon, hartnick$\}$@math.ethz.ch}

\end{document}